\documentclass{article}
\usepackage{amssymb}
\usepackage{amsfonts}
\usepackage{amsmath}

\begin{document}

\begin{center}
{\Large Some properties of Fibonacci numbers, Fibonacci octonions and generalized Fibonacci-Lucas octonions}

\begin{equation*}
\end{equation*}%
Diana SAVIN 
\begin{equation*}
\end{equation*}
\end{center}

\textbf{Abstract. }{\small In this paper we determine some properties of Fibonacci octonions.
Also, we introduce the generalized Fibonacci-Lucas octonions and we investigate
some properties of these elements.}


\medskip


\begin{equation*}
\end{equation*}

\bigskip

\textbf{1. Preliminaries}%
\begin{equation*}
\end{equation*}%
Let $\left( f_{n}\right) _{n\geq 0}$ be the Fibonacci sequence: 
\begin{equation*}
f_{0}=0;f_{1}=1;f_{n}=f_{n-1}+f_{n-2},\;n\geq 2
\end{equation*}%
and let $\left( l_{n}\right) _{n\geq 0}$ be the Lucas sequence: 
\begin{equation*}
l_{0}=2;f_{1}=1;l_{n}=l_{n-1}+l_{n-2},\;n\geq 2
\end{equation*}%
Let $\left( h_{n}\right) _{n\geq 0}$ be the generalized Fibonacci sequence: 
\begin{equation*}
h_{0}=p,h_{1}=q,h_{n}=h_{n-1}+h_{n-2},\;n\geq 2,
\end{equation*}%
where $p$ and $q$ are arbitrary integer numbers. The generalized Fibonacci numbers were introduced of A. F. Horadam in his paper [Ho; 61].\newline
Later A. F. Horadam introduced the Fibonacci quaternions
and generalized Fibonacci quaternions (in the paper [Ho; 63]). In their work [Fl, Sh; 12], C. Flaut and V. Shpakivskyi and later in the paper [Ak, Ko, To; 14], M. Akyigit, H. H
Kosal, M. Tosun gave some properties of the generalized Fibonacci
quaternions. 
In the papers [Fl, Sa; 14] and [Fl, Sa, Io; 13], the authors
introduced the Fibonacci symbol elements and Lucas symbol elements. Moreover, they proved that all these elements determine $\mathbb{Z}$-module structures. 
In the paper [Ke, Ak; 15] O. Kecilioglu, I. Akkus introduced the Fibonacci 
and Lucas octonions and they gave some identities and properties of them.\newline
Quaternion algebras, symbol algebras and octonion algebras have many properties and many applications, as readers can see in [Lam; 04], [Sa; 14], [Sa, Fl, Ci; 09], [Fl, Sa; 15 ], [Fl, Sa; 15 (a)]), [Ke, Ak; 15], [Ke, Ak; 15 (a)], [Ta; 13].
In the paper [Ke, Ak; 15 (a)], O. Kecilioglu, I. Akkus gave some properties of the split Fibonacci 
and Lucas octonions in the octonion algebra $\mathcal {O}$$\left(1,1,-1\right).$ \newline
In this paper we study the Fibonacci octonions in certain generalized octonion algebras.\newline
In the paper [Fl, Sa; 15 (a)], we introduced the generalized Fibonacci - Lucas
quaternions and we determined some properties of these elements. In this paper we introduce the generalized Fibonacci - Lucas octonions and we 
prove that these elements have similar properties with the properties of the generalized Fibonacci - Lucas quaternions.\newline

\begin{equation*}
\end{equation*}

\textbf{2. Properties of the Fibonacci and Lucas numbers}

\begin{equation*}
\end{equation*}

The following properties of Fibonacci and Lucas numbers are known:

\medskip

\textbf{Proposition 2.1.} ([Fib.]) \textit{Let} $(f_{n})_{n\geq 0}$ \textit{be the
Fibonacci sequence} \textit{and let } $(l_{n})_{n\geq 0}$ \textit{be the
Lucas sequence.} \textit{Therefore the following properties hold:}\newline
i) 
\begin{equation*}
f_{n}+f_{n+2}=l_{n+1},\forall ~n\in \mathbb{N};
\end{equation*}%
ii) 
\begin{equation*}
l_{n}+l_{n+2}=5f_{n+1},\forall ~n\in \mathbb{N};
\end{equation*}%
iii) 
\begin{equation*}
f_{n}^{2}+f_{n+1}^{2}=f_{2n+1},\forall ~n\in \mathbb{N};
\end{equation*}%
iv) 
\begin{equation*}
l_{n}^{2}+l_{n+1}^{2}=l_{2n}+l_{2n+2}=5f_{2n+1},\forall ~n\in \mathbb{N};
\end{equation*}%
v) 
\begin{equation*}
l_{n}^{2}=l_{2n}+2\left( -1\right) ^{n},\forall ~n\in \mathbb{N}^{\ast };
\end{equation*}%
vi) 
\begin{equation*}
l_{2n}=5f_{n}^{2}+2\left( -1\right) ^{n},\forall ~n\in \mathbb{N}^{\ast };
\end{equation*}%
vii)
\begin{equation*}
l_{n}+f_{n}=2f_{n+1}.
\end{equation*}
\medskip\newline

\textbf{Proposition 2.2.} ([Fl, Sa; 14], [Fl, Sa, Io; 13]) \textit{Let} $(f_{n})_{n\geq 0}$ \textit{be the Fibonacci sequence} \textit{and let } $(l_{n})_{n\geq 0}$ \textit{be the
Lucas sequence.} \textit{Then:}\newline
i) 
\begin{equation*}
f_{n}+f_{n+3}=2f_{n+2},\forall ~n\in \mathbb{N};
\end{equation*}%
ii) 
\begin{equation*}
f_{n}+f_{n+4}=3f_{n+2},\forall ~n\in \mathbb{N};
\end{equation*}%
ii) 
\begin{equation*}
f_{n}+f_{n+6}=2l_{n+3},\forall ~n\in \mathbb{N};
\end{equation*}%
iv) 
\begin{equation*}
f_{n+4}-f_{n}=l_{n+2},\forall ~n\in \mathbb{N}.
\end{equation*}%
\smallskip
In the following proposition, we will give other properties of the
Fibonacci and Lucas numbers, which will be necessary in the
next proofs.\bigskip\newline
\textbf{Proposition 2.3.} \textit{Let} $(f_{n})_{n\geq 0}$ \textit{be the
Fibonacci sequence} \textit{and} $(l_{n})_{n\geq 0}$ \textit{be the Lucas
sequence} \textit{Then:}\newline
i)
\begin{equation*}
l_{n+4}+l_{n}=3l_{n+2},\forall ~n\in \mathbb{N}.
\end{equation*}%
ii)
\begin{equation*}
l_{n+4}-l_{n}=5f_{n+2},\forall ~n\in \mathbb{N}.
\end{equation*}%
iii)
\begin{equation*}
f_{n}+f_{n+8}=7f_{n+4},\forall ~n\in \mathbb{N}.
\end{equation*}%
\textbf{Proof.} 
i) Using Proposition 2.1 (i) we have:
$$l_{n+4}+l_{n}=f_{n+3}+f_{n+5}+f_{n-1}+f_{n+1}.$$
From Proposition 2.2 (ii) and Proposition 2.1 (i), we obtain:
$$l_{n+4}+l_{n}=3f_{n+1}+3f_{n+3}=3l_{n+2}.$$
ii) Applying Proposition 2.1 (ii), we have:
$$l_{n+4}-l_{n}=\left(l_{n+4}+l_{n+2}\right)-\left(l_{n+2}+l_{n}\right)=$$
$$=5f_{n+3}-5f_{n+1}=5f_{n+2}.$$
iii) $$f_{n}+f_{n+8}=\left(f_{n}+f_{n+4}\right)+\left(f_{n+8}-f_{n+4}\right).$$
Using Proposition 2.2 (ii,iv), we have:
$$f_{n}+f_{n+8}=3f_{n+2}+l_{n+6}.$$
From Proposition 2.1 (i) and Fibonacci recurrence, we obtain:
$$f_{n}+f_{n+8}=3f_{n+2}+f_{n+5}+f_{n+7}=3f_{n+2}+2f_{n+5}+f_{n+6}=3f_{n+2}+3f_{n+5}+f_{n+4}.$$
Using Proposition 2.2 (i), we obtain:
$$f_{n}+f_{n+8}=6f_{n+4}+f_{n+4}=7f_{n+4}.$$

$\square \medskip $ 
\bigskip 
\medskip 
\begin{equation*}
\end{equation*}%
\textbf{3. Fibonacci octonions}%
\begin{equation*}
\end{equation*}
Let $\mathcal {O}_{\mathbb{R}}$$\left(\alpha,\beta,\gamma\right)$ be the generalized octonion algebra over $\mathbb{R}$ with basis $\left\{1, e_{1}, e_{2},..., e_{7}\right\}.$ 
It is known that this algebra is an eight-dimensional non-commutative and non-associative algebra.\\
The multiplication table for the basis of $\mathcal {O}_{\mathbb{R}}$$\left(\alpha,\beta,\gamma\right)$ is
{\footnotesize 
\[
\begin{tabular}{l|llllllll}
$\cdot $ & $1$ & $\,\,e_{1}$ & $\,\,e_{2}$ & $ \,\,\,e_{3}$ & $\,\,e_{4}$ & $\,\,\,e_{5}$ & $\,\,e_{6}$& $\,\,e_{7}$ \\ \hline
 $1$ & $1$ & $\,\,e_{1}$ & $\,\,e_{2}$ & $\,\,\,e_{3}$ & $\,\,e_{4}$ & $\,\,e_{5}$ & $\,\,e_{6}$& $\,\,e_{7}$\\ 
$e_{1}$ & $e_{1}$ & $-\,\,\alpha$ & $\,\,e_{3}$ & $\,\,-\alpha e_{2}$ & $\,\,e_{5}$ & $-\,\,\alpha e_{4}$ & $\,\,-e_{7}$& $\,\,\alpha e_{6}$\\ 
$e_{2}$ & $e_{2}$ & $-e_{3}$ & $-\beta $ & $\,\,\,\beta e_{1}$ & $e_{6} $ & $\,\,\, e_{7}$ & $-\,\,\,\beta e_{4}$ & $-\,\,\,\beta e_{5}$\\ 
$e_{3}$ & $e_{3}$ & $\alpha e_{2}$ & $-\beta e_{1}$ & $-\alpha \beta $ & $e_{7}$ & $-\alpha e_{6}$ & $\beta e_{5}$ & $-\alpha\beta e_{4}$\\
$e_{4}$ & $e_{4}$ & $-e_{5}$ & $-e_{6}$ & $-e_{7}$ & $-\gamma$ & $\gamma e_{1}$ & $\gamma e_{2}$ & $\gamma e_{3}$\\
$e_{5}$ & $e_{5}$ & $\alpha e_{4}$ & $-e_{7}$ & $\alpha e_{6}$ & $-\gamma e_{1}$ & $-\alpha\gamma $ & $-\gamma e_{3}$ & $\alpha\gamma e_{2}$\\
$e_{6}$ & $e_{6}$ & $e_{7}$ & $\beta e_{4}$ & $-\beta e_{5}$ & $-\gamma e_{2}$ & $\gamma e_{3} $ & $-\beta\gamma $ & $-\beta\gamma e_{1}$\\
$e_{7}$ & $e_{7}$ & $-\alpha e_{6}$ & $\beta e_{5}$ & $\alpha\beta e_{4}$ & $-\gamma e_{3}$ & 
$-\alpha\gamma e_{2} $ & $\beta\gamma e_{1}$ & $-\alpha\beta\gamma$
\end{tabular}
\]
}
\smallskip\\
Let $x$$\in$$\mathcal {O}_{\mathbb{R}}$$\left(\alpha,\beta,\gamma\right),$  $x = x_{0}+x_{1}e_{1}+x_{2}e_{2}+x_{3}e_{3}+x_{4}e_{4}+x_{5}e_{5}+x_{6}e_{6}+x_{7}e_{7}$ and its conjugate $\overline{x} = x_{0}-x_{1}e_{1}-x_{2}e_{2}-x_{3}e_{3}-x_{4}e_{4}-x_{5}e_{5}-x_{6}e_{6}-x_{7}e_{7},$ the norm of $x$ is $n\left(x\right)=x\overline{x}=x^{2}_{0}+\alpha x^{2}_{1}+\beta x^{2}_{2}+\alpha\beta x^{2}_{3}+\gamma x^{2}_{4}+\alpha\gamma x^{2}_{5}+ \beta\gamma x^{2}_{6}+\alpha\beta\gamma x^{2}_{7}$$\in$$\mathbb{R.}$\newline
If, for $x$$\in$$\mathcal {O}_{\mathbb{R}}$$\left(\alpha,\beta,\gamma\right),$ we have $n\left(x\right)=0$ if and only if $x = 0,$ then the octonion algebra $\mathcal {O}_{\mathbb{R}}$$\left(\alpha,\beta,\gamma\right),$ is
called a division algebra. Otherwise $\mathcal {O}_{\mathbb{R}}$$\left(\alpha,\beta,\gamma\right)$ is called a split algebra. \newline
Let $K$ be an algebraic number field. It is known the following criterion to decide if an octonion algebra is a division algebra.\newline
\smallskip\\
\textbf{Proposition 3.1.} ([Fl, St; 09]) \textit{A generalized octonion algebra} $\mathcal {O}_{K}$$\left(\alpha,\beta,\gamma\right)$ \textit{is a division algebra if and only if the quaternion algebra} $\mathbb{H}_{K}\left(\alpha,\beta\right)$ \textit{is a division algebra and the equation} $n\left(x\right)=-\gamma$ \textit{does not have solutions in the quaternion algebra} $\mathbb{H}_{K}\left(\alpha,\beta\right).$\newline
\smallskip\\
It is known that the octonion algebra $\mathcal {O}_{\mathbb{R}}$$\left(1,1,1\right)$ is a division algebra and the octonion algebra $\mathcal {O}_{\mathbb{R}}$$\left(1,1,-1\right)$ is a split algebra (see [Ke, Ak; 15 (a)], [Fl, Sh; 15]). In [Fl, Sh; 15] appear the following result, which we allow us to decide if an octonion algebra over $\mathbb{R},$ $\mathcal {O}_{\mathbb{R}}$$\left(\alpha,\beta,\gamma\right)$ is a division algebra or a split algebra.\smallskip\\
\textbf{Proposition 3.2.} ([Fl, Sh; 15]) \textit{We consider the generalized octonion algebra} $\mathcal {O}_{\mathbb{R}}$$\left(\alpha,\beta,\gamma\right),$ \textit{with} $\alpha,\beta,\gamma$$\in$$\mathbb{R}^{*}.$ \textit{Then, there are the following isomorphisms}:\\
i) \textit{if} $\alpha,\beta,\gamma > 0,$ \textit{then the octonion algebra} $\mathcal {O}_{\mathbb{R}}$$\left(\alpha,\beta,\gamma\right)$ \textit{is isomorphic to octonion algebra} $\mathcal {O}_{\mathbb{R}}$$\left(1,1,1\right);$\\
ii) \textit{if} $\alpha,\beta > 0, \gamma <0$ or $\alpha,\gamma > 0, \beta <0$ or $\alpha <0,\beta,\gamma > 0$ or $\alpha >0,\beta,\gamma < 0$ or $\alpha,\gamma < 0, \beta >0$ or $\alpha,\beta < 0, \gamma >0$ or $\alpha,\beta,\gamma <0$ \textit{then the octonion algebra} $\mathcal {O}_{\mathbb{R}}$$\left(\alpha,\beta,\gamma\right)$ \textit{is isomorphic to the octonion algebra} $\mathcal {O}_{\mathbb{R}}$$\left(1,1,-1\right).$\newline
\smallskip\\
Let $n$ be an integer, $n\geq0.$ In the paper [Ke, Ak; 15], O. Kecilioglu, I. Akkus, introduced the Fibonacci octonions:
$$F_n=f_{n}+ f_{n+1}e_{1}+ f_{n+2}e_{2}+f_{n+3}e_{3}+ f_{n+4}e_{4}+ f_{n+5}e_{5}+ f_{n+6}e_{6}+ f_{n+7}e_{7},$$ 
where $f_{n}$ is $n^{th}$ Fibonacci number.\newline
\smallskip\\
Now, we consider the generalized octonion algebra $\mathcal {O}_{\mathbb{R}}$$\left(\alpha,\beta,\gamma\right),$ with $\alpha,\beta,\gamma$ in arithmetic progression,
$\alpha=a+1,$ $\beta=2a+1,$ $\gamma=3a+1,$ where $a$$\in$$\mathbb{R}.$ \\
In the following, we calculate the norm of a Fibonacci octonion in this octonion algebra.\newline
\smallskip\\
\textbf{Proposition 3.3.} \textit{Let} $a$ \textit{be a real number and let}  $F_{n}$ \textit{be the}
$n$-\textit{th generalized Fibonacci octonion. Then the norm of} $F_{n}$ \textit{in the generalized octonion algebra} $\mathcal {O}_{\mathbb{R}}$$\left(a+1, 2a+1, 3a+1,\right)$ \textit{is}:\\
$$n\left(F_{n}\right)=f_{2n+6}\left(79a^{2}+ 46a+ \frac{174a^{3}-4a}{5}\right)+$$
$$+f_{2n+7}\left(130a^{2}+ 84a + 21+ \frac{282a^{3}+ 8a}{5}\right)+ \left(-1\right)^{n}\left(
4a^{2} + \frac{12a^{3}+ 8a}{5}\right).$$
\textbf{Proof.}
$$ n\left(F_{n}\right)=f_{n}^{2}+\left(a+1\right)f_{n+1}^{2}+\left(2a+1\right)f_{n+2}^{2}+
\left(a+1\right)\left(2a+1\right)f_{n+3}^{2}+$$

$$+\left(3a+1\right)f_{n+4}^{2}
+\left(a+1\right)\left(3a+1\right)f_{n+5}^{2}+$$

$$+\left(2a+1\right)\left(3a+1\right)f_{n+6}^{2}
+\left(a+1\right)\left(2a+1\right)\left(3a+1\right)f_{n+7}^{2}=$$

$$=f_{n}^{2}+f_{n+1}^{2}+f_{n+2}^{2}+f_{n+3}^{2}+f_{n+4}^{2}+f_{n+5}^{2}+f_{n+6}^{2}+
f_{n+7}^{2}+$$

$$+ a\left(f_{n+1}^{2}+2f_{n+2}^{2}+3f_{n+3}^{2}+3f_{n+4}^{2}+4f_{n+5}^{2}+5f_{n+6}^{2}+
6f_{n+7}^{2}\right)+$$
$$+ a^{2}\left(2f_{n+3}^{2}+3f_{n+5}^{2}+6f_{n+6}^{2}+
11f_{n+7}^{2}\right)+ 6a^{3}f_{n+7}^{2}=$$
$$=S_{1}+S_{2}+S_{3}+ 6a^{3}f_{n+7}^{2}, \eqno(3.1)$$
where, we denoted with $S_{1}=f_{n}^{2}+f_{n+1}^{2}+f_{n+2}^{2}+f_{n+3}^{2}+f_{n+4}^{2}+f_{n+5}^{2}+f_{n+6}^{2}+f_{n+7}^{2},$ $S_{2}=a\left(f_{n+1}^{2}+2f_{n+2}^{2}+3f_{n+3}^{2}+3f_{n+4}^{2}+4f_{n+5}^{2}+5f_{n+6}^{2}+
6f_{n+7}^{2}\right),$ $S_{3}=a^{2}\left(2f_{n+3}^{2}+3f_{n+5}^{2}+6f_{n+6}^{2}+
11f_{n+7}^{2}\right).$\newline
Now, we calculate $S_{1}, S_{2}, S_{3}.$\\
Using [Ke, Ak; 15] (p.3), we have 
$$ S_{1}=f_{8}f_{2n+7}=21f_{2n+7.} \eqno(3.2)$$
Applying Proposition 2.1 (iii) and Proposition 2.1 (i), we
have:

$$S_{2}=a\left(f_{n+1}^{2}+2f_{n+2}^{2}+3f_{n+3}^{2}+3f_{n+4}^{2}+4f_{n+5}^{2}+5f_{n+6}^{2}+
6f_{n+7}^{2}\right)=$$ 

$S_{2}=6a\left(f_{n+6}^{2}+f_{n+7}^{2}\right) -a\left(f_{n+6}^{2}+f_{n+5}^{2}\right) + 5a\left(f_{n+5}^{2}+f_{n+4}^{2}\right) - 2a\left(f_{n+4}^{2}+f_{n+3}^{2}\right)+$\\
$$+  4af_{n+3}^{2}+ a\left(f_{n+2}^{2}+f_{n+3}^{2}\right) + a\left(f_{n+1}^{2}+f_{n+2}^{2}\right)=$$
$$=a\cdot \left(6f_{2n+13}-f_{2n+11} +5f_{2n+9}- 2f_{2n+7}+ 4f_{n+3}^{2}+f_{2n+5}+f_{2n+3}\right)=$$
$$=a\cdot \left[6f_{2n+13}-f_{2n+11} + 7f_{2n+9}- 2\left(f_{2n+7}+f_{2n+9}\right)+ 4f_{n+3}^{2}+ l_{2n+4}\right] = $$
$$=a\cdot \left[6\left(f_{2n+9}+f_{2n+13}\right)+\left(f_{2n+9}-f_{2n+11}\right)+ l_{2n+4} -
2l_{2n+8} + 4\cdot\frac{l_{2n+6}-2\left(-1\right)^{n+3}}{5}\right].$$
From Proposition 2.2 (ii), Proposition 2.3 (ii) and Proposition 2.1 (i), we have:
$$S_{2}=a\cdot \left[18f_{2n+11}- f_{2n+10}-\left(l_{2n+8}-l_{2n+4}\right)-l_{2n+8}+
4\cdot\frac{l_{2n+6}+2\left(-1\right)^{n}}{5}\right]=$$ 
$$=a\cdot \left[18f_{2n+11}- f_{2n+10}-5f_{2n+6}-f_{2n+7}-f_{2n+9}+ 4\cdot\frac{l_{2n+6}+2\left(-1\right)^{n}}{5}\right].$$
Using several times the recurrence of Fibonacci sequence and Proposition 2.1 (vii), we obtain:
$$S_{2}=a\cdot \left[46f_{2n+6}+84f_{2n+7}+ 4\cdot\frac{2f_{2n+7}-f_{2n+6}+2\left(-1\right)^{n}}{5}\right].\eqno(3.3)$$
Applying Proposition 2.1 (iii) and Proposition 2.1 (vi,i), we have:
$$S_{3}=a^{2}\cdot\left(2f_{n+3}^{2}+3f_{n+5}^{2}+6f_{n+6}^{2}+
11f_{n+7}^{2}\right)=$$
$$=a^{2}\cdot \left[2\left(f_{n+3}^{2}+f_{n+4}^{2}\right)-2\left(f_{n+4}^{2}+f_{n+5}^{2}\right) + 5\left(f_{n+5}^{2}+f_{n+6}^{2}\right) + \left(f_{n+6}^{2}+f_{n+7}^{2}\right) + 10f_{n+7}^{2}\right]=$$
$$=a^{2}\cdot \left[2f_{2n+7}-2f_{2n+9} + 5f_{2n+11}+f_{2n+13}+2l_{2n+14}-4\left(-1\right)^{n+7}\right]=$$
$$=a^{2}\cdot \left[-2f_{2n+8}+5f_{2n+11}+3f_{2n+13}+2f_{2n+15}+4\left(-1\right)^{n}\right].$$
From Proposition 2.3 (iii) and the recurrence of Fibonacci sequence, we have:
$$S_{3}=a^{2}\cdot \left[-2f_{2n+7}-2f_{2n+6}+5f_{2n+11}+3f_{2n+13}+14f_{2n+11}-2f_{2n+7}+4\left(-1\right)^{n}\right]=$$
$$=a^{2}\cdot \left[-4f_{2n+7}-2f_{2n+6}+19f_{2n+11}+3f_{2n+13}+4\left(-1\right)^{n}\right]=$$
$$=a^{2}\cdot \left[-4f_{2n+7}-2f_{2n+6}+134f_{2n+7}+81f_{2n+6}+4\left(-1\right)^{n}\right].$$
Therefore, we obtained that:
$$S_{3}=a^{2}\cdot \left[79f_{2n+6}+130f_{2n+7}+4\left(-1\right)^{n}\right].\eqno(3.4)$$
From Proposition 2.1 (vi,i), we have:
$$6a^{3}f_{n+7}^{2}=\frac{6a^{3}}{5}\cdot \left[l_{2n+14}-2\cdot\left(-1\right)^{n+7}\right]=
\frac{6a^{3}}{5}\cdot \left[f_{2n+13}+f_{2n+15}+2\cdot\left(-1\right)^{n}\right].$$
Applying Proposition 2.3 (iii) and the recurrence of Fibonacci sequence many times, we
have:
$$6a^{3}f_{n+7}^{2}=\frac{6a^{3}}{5}\cdot \left[7f_{2n+9}-f_{2n+5}+7f_{2n+11}-f_{2n+7}+2\cdot\left(-1\right)^{n}\right]=$$
$$=\frac{6a^{3}}{5}\cdot \left[29f_{2n+6}+ 47f_{2n+7}+2\cdot\left(-1\right)^{n}\right].\eqno(3.5)$$
From the relations (3.1), (3.2), (3.3), (3.4), (3.5), we have:
$$n\left(F_{n}\right)=21f_{2n+7}+a\cdot \left[46f_{2n+6}+84f_{2n+7}+ 4\cdot\frac{2f_{2n+7}-f_{2n+6}+2\left(-1\right)^{n}}{5}\right]+$$
$$+a^{2}\cdot \left[79f_{2n+6}+130f_{2n+7}+4\left(-1\right)^{n}\right]
+\frac{6a^{3}}{5}\cdot \left[29f_{2n+6}+ 47f_{2n+7}+2\cdot\left(-1\right)^{n}\right].$$
Therefore, we get
:
$$n\left(F_{n}\right)=f_{2n+6}\left(79a^{2}+ 46a+ \frac{174a^{3}-4a}{5}\right)+$$
$$+f_{2n+7}\left(130a^{2}+ 84a + 21+ \frac{282a^{3}+ 8a}{5}\right)+ \left(-1\right)^{n}\left(
4a^{2} + \frac{12a^{3}+ 8a}{5}\right)$$
$\square \medskip $\newline
 We obtain immediately the following remark:\newline
\smallskip\\
\textbf{Remark 3.1.} \textit{If} $a$ \textit{is a real number,} $a< -1,$ \textit{then, the generalized octonion algebra} $\mathcal {O}_{\mathbb{R}}$$\left(a+1, 2a+1, 3a+1\right)$ \textit{is a split algebra.}\newline
\smallskip\\
\textbf{Proof.} Using Proposition 3.2 (ii) and the fact that the octonion algebra $\mathcal {O}_{\mathbb{R}}$$\left(1,1,-1\right)$ is a split algebra, it results that, if $a< -1,$ the generalized octonion algebra $\mathcal {O}_{\mathbb{R}}$$\left(a+1, 2a+1, 3a+1\right)$ is a split algebra.
\newline
\smallskip\\
For example, for $a=-4$ we obtain the generalized octonion algebra $\mathcal {O}_{\mathbb{R}}$$\left(-3, -7, -11 \right).$ From Remark 3.1 it results that this is a split algebra (another way for to  prove that this algebra is a split algebra is to remark that the equation $n\left(x\right)=11$ has solutions in the quaternion algebra $\mathbb{H}_{K}\left(-3,-7\right)$ and then we apply Proposition 3.1.\newline 
Now, we want to determine how many Fibonacci octonions invertible are in the octonion algebra $\mathcal {O}_{\mathbb{R}}$$\left(-3, -7, -11 \right).$ Applying Proposition 3.3, we obtain that $n\left(F_{n}\right)=-1144f_{2n+6}-1851f_{2n+7}-96\left(-1\right)^{n},$ $n\in$$\mathbb{N}.$ Using that $f_{2n+6}, f_{2n+7} > 0,$ $\left(\forall\right) n$$\in$$\mathbb{N},$ it results that $n\left(F_{n}\right)<0,$ $\left(\forall\right) n$$\in$$\mathbb{N},$ therefore, in the split octonion algebra $\mathcal {O}_{\mathbb{R}}$$\left(-3, -7, -11 \right)$ all Fibonacci octonions are  invertible. \newline
For $a = -2,$ after a few calculations, we also get that in the split octonion algebra $\mathcal {O}_{\mathbb{R}}$$\left(-1, -3, -5 \right)$ all Fibonacci octonions are  invertible. \newline
From the above, the following question arises: how many invertible Fibonacci octonions there are in the octonion algebra
$\mathcal {O}_{\mathbb{R}}$$\left(a+1, 2a+1, 3a+1\right),$ with $a<-1?$ We get the following result:\newline
\smallskip\\
\textbf{Proposition 3.4.} \textit{Let} $a$ \textit{be a real number,} $a\leq -2$ \textit{and let} $\mathcal {O}_{\mathbb{R}}$$\left(a+1, 2a+1, 3a+1\right)$ \textit{be a generalized octonion algebra. Then, in this algebra, all Fibonacci octonions are  invertible elements.}\newline
\smallskip\\
\textbf{Proof.} It is sufficient to prove that $n\left(F_{n}\right)\neq 0,$ $\left(\forall\right) n$$\in$$\mathbb{N}.$
Using Proposition 3.3., we have:
$$n\left(F_{n}\right)=f_{2n+6}\left(79a^{2}+ 46a+ \frac{174a^{3}-4a}{5}\right)+$$
$$+f_{2n+7}\left(130a^{2}+ 84a + 21+ \frac{282a^{3}+ 8a}{5}\right)+ \left(-1\right)^{n}\left(
4a^{2} + \frac{12a^{3}+ 8a}{5}\right) \Leftrightarrow$$
$$n\left(F_{n}\right)=f_{2n+6}\cdot\frac{174a^{3}+395a^{2}+226a}{5}+$$
$$+f_{2n+7}\cdot\frac{282a^{3}+650a^{2}+428a+105}{5}
+\left(-1\right)^{n}\cdot
 \frac{12a^{3}+20a^{2}+ 8a}{5}.$$
After a few calculations, we obtain:
$$n\left(F_{n}\right)=f_{2n+6}\cdot\frac{a\left(a+2\right)\left(174a+47\right)+132a}{5}+$$
$$+f_{2n+7}\cdot\frac{2\left(a+2\right)\left(141a^{2}+43a+126\right)-407}{5}+
\left(-1\right)^{n}\cdot\frac{4a\left(a+2\right)\left(3a-1\right)+16a}{5}.$$
We remark that $141a^{2}+43a+126 >0$ $\left(\forall\right) a$$\in$$\mathbb{R}$ (since $\Delta < 0$)and \newline
$\frac{a\left(a+2\right)\left(174a+47\right)+132a}{5}<0,$ $\left(\forall\right)$$a\leq -2,$ $\frac{2\left(a+2\right)\left(141a^{2}+43a+126\right)-407}{5}<0,$ $\left(\forall\right)$$a\leq -2,$ $\frac{4a\left(a+2\right)\left(3a-1\right)+16a}{5}<0,$ $\left(\forall\right)$$a\leq -2.$ 
Since $f_{2n+6}, f_{2n+7}>0$ $\left(\forall\right) n$$\in$$\mathbb{N},$ we obtain that
$n\left(F_{n}\right)<0,$ $\left(\forall\right)$$a\leq -2,$ $n$$\in$$\mathbb{N}$ (even if $n$ is an odd number). This implies that, in the generalized octonion algebra $\mathcal {O}_{\mathbb{R}}$$\left(a+1, 2a+1, 3a+1,\right),$ with $a\leq -2,$ all Fibonacci octonions are  invertible.
$\square \medskip $\newline
Now, we wonder what happens with the Fibonacci octonions in the generalized octonion algebra $\mathcal {O}_{\mathbb{R}}$$\left(a+1, 2a+1, 3a+1\right),$ when $a$$\in$$\left(-2,-1\right)$?
All Fibonacci octonions in a such octonion algebra are  invertible or there are Fibonacci octonions zero divisors?\newline
For example, for $a=-\frac{3}{2},$ using Proposition 3.3, we have that the norm of a Fibonacci octonion in the octonion algebra $\mathcal {O}_{\mathbb{R}}$$\left(-\frac{1}{2}, -2, -\frac{7}{2}\right)$ is $n\left(F_{n}\right)=-\frac{7}{2}f_{2n+6} + \frac{831}{2}f_{2n+7}   -\frac{3}{2}\cdot\left(-1\right)^{n}$$>$$0,$ $\left(\forall\right) n$$\in$$\mathbb{N}^{*}.$ This implies that in the generalized octonion algebra $\mathcal {O}_{\mathbb{R}}$$\left(-\frac{1}{2}, -2, -\frac{7}{2}\right)$ all Fibonacci octonions are  invertible.\newline
In the future, we will study if this fact is true in each generalized octonion algebra $\mathcal {O}_{\mathbb{R}}$$\left(a+1, 2a+1, 3a+1\right)$, with $a$$\in$$\left(-2,-1\right).$

\begin{equation*}
\end{equation*}%
\textbf{4. Generalized Fibonacci-Lucas octonions}%
\begin{equation*}
\end{equation*}
In the paper [Fl, Sa; 15 (a)], we introduced the generalized Fibonacci - Lucas
numbers, namely:
if $n$ is an arbitrary positive integer and $p,q$ be two arbitrary
integers, the sequence $\left(g_{n}\right)_{n\geq 1},$ where 
\begin{equation*}
g_{n+1}=pf_{n}+ql_{n+1},\;n\geq 0
\end{equation*}%
is called  \textit{the sequence of the generalized Fibonacci-Lucas numbers}. For not make  confusions, we will  use the notation $g_{n}^{p,q}$ instead of $g_{n}.$ \newline 
Let $\mathcal {O}_{\mathbb{Q}}$$\left(\alpha,\beta,\gamma\right)$ be the generalized octonion algebra over $\mathbb{Q}$ with the basis $\left\{1, e_{1}, e_{2},..., e_{7}\right\}.$ 
We define the $n$-%
\textit{th generalized Fibonacci-Lucas octonion} to be the element of the
form 
\begin{equation*}
G_{n}^{p,q}=g_{n}^{p,q}\cdot 1+g_{n+1}^{p,q}\cdot e_{1}+g_{n+2}^{p,q}\cdot
e_{2}+g_{n+3}^{p,q}\cdot e_{3}+g_{n+4}^{p,q}\cdot e_{4}+g_{n+5}^{p,q}\cdot e_{5}+g_{n+6}^{p,q}\cdot e_{6}+g_{n+7}^{p,q}\cdot e_{7}.
\end{equation*}%
We wonder what algebraic structure determine the generalized Fibonacci-Lucas octonions. First, we make the following remark.
\newline
\smallskip\\
\textbf{Remark 4.1.} \textit{Let} $n, p, q$ \textit{three arbitrary positive integers,} $p,q\geq0.$ \textit{Then, the} $n$-%
\textit{th generalized Fibonacci-Lucas octonion} $G_{n}^{p,q}=0$ \textit{if and only if} $p=q=0.$\newline
\smallskip\\
\textbf{Proof.} $"\Rightarrow"$ If $G_{n}^{p,q}=0,$ it results $g_{n}^{p,q}=g_{n+1}^{p,q}=...=g_{n+7}^{p,q}=0.$ This implies that $g_{n-1}^{p,q}=...=g_{1}^{p,q}=g_{1}^{p,q}=0.$ We obtain immediately that $q=0$ and $p=0.$\newline
$"\Leftarrow"$ is trivial.\newline
\smallskip\\
In the paper [Fl, Sa; 15 (a)], we proved the following properties of the generalized Fibonacci - Lucas numbers:\newline
\smallskip\\
\textbf{Remark 4.2.} \textit{Let} $n,m\in \mathbb{N}^{\ast },$ $a,b,p,q,p^{^{\prime }},q^{^{\prime }}\in \mathbb{Z}.$ \textit{Then, we have:}\newline
i) \begin{equation*}
ag_{n}^{p,q}+bg_{m}^{p^{^{\prime }},q^{^{\prime
}}}=g_{n}^{ap,aq}+g_{m}^{bp^{^{\prime }},bq^{^{\prime }}};
\end{equation*}%
ii) \begin{equation*}
5g_{n}^{p,q}\cdot 5g_{m}^{p^{^{\prime }},q^{^{\prime }}}=
5g_{m+n-2}^{5p^{^{\prime }}q,pp^{^{\prime }}}+5g_{m+n-1}^{5p^{^{\prime
}}q,0}+
\end{equation*}%
\begin{equation*}
5g_{n-m}^{5p^{^{\prime }}q\cdot \left( -1\right) ^{m},pp^{^{\prime
}}\cdot \left( -1\right) ^{m}}+5g_{n-m+1}^{5p^{^{\prime }}q\cdot \left(
-1\right) ^{m},0}+
5g_{m+n}^{5pq^{^{\prime }},5qq^{^{\prime }}}+5g_{n-m}^{5pq^{^{\prime
}}\cdot \left( -1\right) ^{m},5qq^{^{\prime }}\cdot \left( -1\right) ^{m}}.
\end{equation*}%
Using these remark we can prove the followings:\newline
\smallskip\\
\textbf{Theorem 4.1.} \textit{Let} $A$ \textit{and} $B$ \textit{be} \textit{the sets} 
\begin{equation*}
A=\left\{ \sum\limits_{i=1}^{n}G_{n_{i}}^{p_{i},q_{i}}|n\in \mathbb{N}%
^{\ast }, p_{i},q_{i}\in \mathbb{Z},(\forall )i=\overline{1,n}\right\},
\end{equation*}%
\begin{equation*}
B=\left\{ \sum\limits_{i=1}^{n}5G_{n_{i}}^{p_{i},q_{i}}|n\in \mathbb{N}%
^{\ast },p_{i},q_{i}\in \mathbb{Q},(\forall )i=\overline{1,n}\right\} \cup
\left\{ 1\right\}.
\end{equation*}%
\textit{Then, the following statements are true:}\newline
\textit{i)} $A$ \textit{is a free} $\mathbb{Z}$- \textit{submodule of rank} $8$ \textit{of the generalized octonions algebra} $\mathcal {O}_{\mathbb{Q}}$$\left(\alpha,\beta,\gamma\right);$\newline  
\textit{ii)} $B$ \textit{with octonions addition and multiplication, is a unitary non-associative subalgebra of the generalized octonions algebra} $\mathcal {O}_{\mathbb{Q}}$$\left(\alpha,\beta,\gamma\right).$
\newline
\smallskip\\
\textbf{Proof.} (i) Using Remark 4.2, it results immediately that $$aG_{n}^{p,q}+bG_{m}^{p^{^{\prime }},q^{^{\prime
}}}=G_{n}^{ap,aq}+G_{m}^{bp^{^{\prime }},bq^{^{\prime }}},\ \left(\forall\right) m,n\in \mathbb{N}^{*}, a,b,p,q,p^{{\prime }},q^{{\prime }}\in\mathbb{Z}.$$ 
Moreover, applying Remark 4.1, it results that $0$$\in$$A.$\newline
These implies that $A$ is a $\mathbb{Z}$ - submodule of the generalized octonions algebra 
$\mathcal {O}_{\mathbb{Q}}$$\left(\alpha,\beta,\gamma\right).$. Since $\left\{1, e_{1}, e_{2},..., e_{7}\right\}$ is a basis of $A$, it results
that $A$ is a free $\mathbb{Z}$-module of rank $8$.\newline
(ii) From Remark 4.2 (ii), it results immediately that $5G_{m}^{p,q}\cdot 5G_{n}^{p^{^{\prime
}},q^{^{\prime
}}}$$\in$$B,$ $\left(\forall\right)$ $m,n$$\in$$\mathbb{N}^{*},$ 
$p,q,p^{{\prime }},q^{{\prime }}$$\in$$\mathbb{Z}.$ Using this fact and a similar reason that in the proof of (i), it results that $B$ is a unitary non-associative subalgebra of the generalized octonions algebra $\mathcal {O}_{\mathbb{Q}}$$\left(\alpha,\beta,\gamma\right).$

$\square \medskip $\newline
\textbf{Acknowledgements.} The author is very grateful her colleague Cristina Flaut,
for helpful discutions about the octonion algebras, which helped the
author to improve this paper.
\newline
\smallskip\\
\textbf{References}
\begin{equation*}
\end{equation*}%
[Ak, Ko, To; 14] M. Akyigit, HH Kosal, M. Tosun, \textit{\ Fibonacci
Generalized Quaternions }, Advances in Applied Clifford Algebras, vol. \textbf{24},
issue: 3 (2014), p. 631-641.\newline
[Fl, St; 09] C. Flaut, M. Stefanescu,
 \textit{Some equations over generalized quaternion and octonion division algebras}, Bull. Math. Soc. Sci. Math. Roumanie, \textbf{52(4)(100)} (2009), p. 427-439. \newline
[Fl, Sa, Io; 13] C. Flaut, D. Savin, G. Iorgulescu, \textit{Some properties
of Fibonacci and Lucas symbol elements}, Journal of Mathematical Sciences:
Advances and Applications, vol. \textbf{20} (2013), p. 37-43 \newline
[Fl, Sa; 14] C. Flaut, D. Savin, \textit{Some properties of the symbol
algebras of degree} $3$, Math. Reports, vol. \textbf{16(66)}(3)(2014),
p.443-463.\newline
[Fl, Sa; 15] C. Flaut, D. Savin, \textit{Some examples of division symbol algebras of degree} $3$ \textit{and} $5$, Carpathian J. Math., vol. 31, No. 2 (2015), p. 197-204.
\newline
[Fl, Sa; 15 (a)] C. Flaut, D. Savin, \textit{Quaternion Algebras and Generalized Fibonacci-Lucas Quaternions}, accepted for publication in Advances in Applied Clifford Algebras
(http://link.springer.com/article/10.1007
).\newline
[Fl, Sh; 12] C. Flaut, V. Shpakivskyi, \textit{On Generalized Fibonacci
Quaternions and Fibonacci-Narayana Quaternions,} Advances in Applied
Clifford Algebras, vol. \textbf{23} issue 3 (2013), p. 673-688.\newline
[Fl, Sh; 15] C. Flaut, V. Shpakivskyi, \textit{, An Efficient Method for Solving Equations in Generalized Quaternion and Octonion Algebras,} accepted for publication in Advances in Applied
Clifford Algebras \newline
(http://link.springer.com/article/10.1007/s00006-014-0493-x).\newline
[Ho; 61] A. F. Horadam, \textit{A Generalized Fibonacci Sequence}, Amer.
Math. Monthly, \textbf{68}(1961), 455-459.\newline
[Ho; 63] A. F. Horadam, \textit{Complex Fibonacci Numbers and Fibonacci
Quaternions}, Amer. Math. Monthly, \textbf{70}(1963), 289-291.\newline
[Ke, Ak; 15] O. Kecilioglu, I. Akkus, \textit{The Fibonacci Octonions},
Adv. Appl. Clifford Algebras, vol. \textbf{25} issue 1 (2015), p. 151-158.\newline
[Ke, Ak; 15 (a)] O. Kecilioglu, I. Akkus, \textit{Split Fibonacci and Lucas Octonions},
accepted for publication in Adv. Appl. Clifford Algebras \newline
(http://link.springer.com/article/10.1007/s00006-014-0515-8).\newline
[Lam; 04] T. Y. Lam, \textit{Introduction to Quadratic Forms over Fields,}
American Mathematical Society, 2004.\newline
[Ol; 14] G. Olteanu, \textit{Baer-Galois connections and applications}, Carpathian J. Math., 30 (2014), No. 2, 225–229. \newline
[Sa, Fl, Ci; 09] D. Savin, C. Flaut, and C. Ciobanu, \textit{Some properties of the symbol algebras}, Carpathian J. Math., 25 (2009), No.2, 239-245. \newline
[Sa; 14] D. Savin, \textit{About some split central simple algebras}, An.
Stiin. Univ. "Ovidius" Constanta, Ser. Mat, \textbf{22} (1) (2014), p.
263-272.\newline
[Ta; 13] M. Tarnauceanu, \textit{A characterization of the quaternion group,} An.
Stiin. Univ. "Ovidius" Constanta, Ser. Mat, \textbf{21} (1) (2013), p. 209-214.\newline
[Fib.] http://www.maths.surrey.ac.uk/hosted-sites/R.Knott/Fibonacci/fib.html%

\begin{equation*}
\end{equation*}

Diana SAVIN

{\small Faculty of Mathematics and Computer Science, Ovidius University,}

{\small Bd. Mamaia 124, 900527, CONSTANTA, ROMANIA }

{\small http://www.univ-ovidius.ro/math/}

{\small e-mail: \ savin.diana@univ-ovidius.ro, \ dianet72@yahoo.com}\bigskip
\bigskip

\end{document}